\documentclass[a4paper,10pt]{article}
\usepackage{amssymb}
\usepackage{latexsym}
\begin{document}
\newtheorem{prop}{Proposition}
\newtheorem{thm}{Theorem}
\newtheorem{rem}{Remark}
\newtheorem{cor}{Corollary}
\newtheorem{lemma}{Lemma}
\newtheorem{definition}{Definition}
\title{Definition and Self-Adjointness of the Stochastic \\ Airy Operator
\thanks{Dedicated to Professor Leonid Pastur on the occasion 
of his 75th birthday} }
\author{Nariyuki {\sc Minami}\footnote{School of medicine, Keio 
University, Hiyoshi 4-1-1, Kohoku-ku, Yokohama 223-8521, Japan. 
e-mail: minami@a5.keio.jp}}
\date{}
\maketitle
\begin{abstract}
In this note, it is shown that the stochastic Airy operator, which is the 
\lq\lq Schr\"odinger operator\rq\rq on the half-line whose potential 
term consists of Gaussian white noise plus a linear term tending to 
$+\infty$, can naturally be defined as a generalized Sturm-Liouville 
operator and that it is self-adjoint and has purely discrete spectrum 
with probability one. Thus \lq\lq 
stochastic Airy spectrum\rq\rq of Ram\'irez, Rider and Vir\'ag 
is the spectrum of an operator in the ordinary sense of the word.


\end{abstract}

\section{Introduction.}
In \cite{DE}, Dumitriu and Edelman considered the following random 
matrix
\begin{equation}
H_n^{\beta}=
\frac1{\sqrt{\beta}}
\left[
\begin{array}{cccccccc} 
g_1 & \chi_{(n-1)\beta} & 0 &   & \ldots &  & \ldots & 0\\ 
\chi_{(n-1)\beta} & g_2 & \chi_{(n-2)\beta} & 0 &   & \ldots &  & 0\\
\ &\ldots &  &\ldots &  &\ldots &  &\ldots \\
\ &\ldots &  &\ldots &  &\ldots &  &\ldots \\
0 & \ldots &  & \ldots &  & \chi_{2\beta} & g_{n-1} & \chi_{\beta} \\
0 & \ldots &  & \ldots &  & 0 & \chi_{\beta} & g_n
\end{array}
\right],
\label{eqn:1}
\end{equation}
where $\beta>0$ is a real number and the random variables 
$g_1,\ldots,g_n$ obey the normal distribution $N(0,2)$, whereas 
$\chi_{j\beta}$ ($j=1,\ldots,n-1$) obeys the $\chi$-distribution 
with parameter $j\beta$, all random variables being independent 
each other. Here by definition, a random variable $X$ obeys the 
$\chi$-distribution if and only if $X^2$ obeys the $\chi^2$-
distribution. They called $H_n^{\beta}$ the $\beta$-ensemble, 
and proved that the $n$ eigenvalues of $H_n^{\beta}$ have the 
joint probability density of the form
\begin{equation}
P_n^{\beta}(\lambda_1,\ldots,\lambda_n)
=\frac1{Z_n^{\beta}}e^{-\frac{\beta}{4}\sum_{k=1}^n\lambda_k^2}
\prod_{j<k}\vert\lambda_j-\lambda_k\vert^{\beta},
\label{eqn:2}
\end{equation}
$Z_n^{\beta}$ being the normalizing constant. For the special values 
$\beta=1,2,4$, (\ref{eqn:2}) represents the joint probability density 
of the eigenvalues of the random matrix GOE, GUE and GSE 
respectively (see e.g. \cite{Me}). It was found out by Trotter 
\cite{T} that these three well studied matrix ensembles can be 
transformed into the tridiagonal form (\ref{eqn:1}), and in this form, 
the general $\beta$-ensemble interpolate and extrapolate those 
three ensembles. 

Concerning the limiting distribution of the largest $k$ eigenvalues 
$\lambda_1^{(n)}\geq\cdots\geq\lambda_k^{(n)}$ of $H_n^{\beta}$, 
Ram\'irez, Rider and Vir\'ag \cite{RRV} proved that as $n\to\infty$, 
the joint distribution of $\{n^{1/6}(2\sqrt{n}-\lambda_j^{(n)})\}_{j=1}^k$ 
converges to that of the smallest $k$ eigenvalues $\Lambda_0\leq
\cdots\leq\Lambda_{k-1}$ of the \lq\lq Schr\"odinger operator\rq\rq 
\begin{equation}
H=-\frac{d^2}{dt^2}+t+\frac{2}{\sqrt{\beta}}B_{\omega}'(t)\ ,\quad 
t\geq0
\label{eqn:3}
\end{equation}
with Dirichlet boundary condition at $t=0$, where 
$\{B_{\omega}(t)\}_{t\geq0}$ is the standard Brownian motion defined 
on a probability space $(\Omega,{\cal F},\mathbf{P})$, and 
$B'_{\omega}(t)$ is the formal derivative of its sample function.

Since sample functions of the Brownian motion are not differentiable, 
the expression (\ref{eqn:3}) itself needs a justification. Ram\'irez et. al. 
consider $H$ to be a random linear mapping from the function 
space
$$H^1_{loc}(\mathbf{R}_+):=
\{f;\ f(\cdot)\ \mbox{is absolutely continuous on}\ [0,\infty)\ 
\mbox{and}\ f'\in L^2_{loc}(\mathbf{R}_+)\}$$
into the space of Schwartz distribution. They further define 
the eigenvalues and eigenfunctions of $H$ to be pairs 
$(\lambda,f)\in\mathbf{R}\times L^{\ast}$ such that $Hf=\lambda f$ 
holds, 
both sides of which being interpreted as Schwartz distributions. 
Here the function space $L^{\ast}$ is defined by 
\begin{eqnarray}
L^{\ast}=\{f&;&\ f(\cdot)\ \mbox{is absolutely continuous on}\ 
[0,\infty),\ f(0)=0,\ \mbox{and} \nonumber\\
& &
\int_0^{\infty}\{\vert f'(t)\vert^2+(1+t)\vert f(t)\vert^2\}dt<\infty\}\ .
\label{eqn:4}
\end{eqnarray}

Actually, Schwartz distribution is irrelevant in formulating the 
eigenvalue problem for $H$, and it turns out that $H$ is realized 
quite naturally as a self-adjoint operator in the Hilbert space 
$L^2(\mathbf{R}_+)$. Namely we have the following

\begin{thm}
For each $\omega$, $H$ can be realized as a closed symmetric 
operator in $L^2(\mathbf{R}_+)$, which is self-adjoint with probability 
one and has purely discrete spectrum.
\end{thm}

Thus what Ram\'irez et.al. called  \lq\lq stochastic Airy 
spectrum\rq\rq is the spectrum of a self-adjoint 
operator in the ordinary sense of the word. We shall prove Theorem 1 
as an immediate corollary of the following more general result.

\begin{thm}
Let $p(t)\geq0$ be a real valued continuous function on $\mathbf{R}_+$ 
such that $\liminf_{t\to\infty}p(t)/t^{\alpha}>0$ for some $\alpha>0$. 
Further let $\{X_{\omega}(t)\}$ be a fractional Brownian motion 
of Hurst parameter $h\in(0,1)$ defined on a probability space 
$(\Omega,{\cal F},\mathbf{P})$. Then under the Dirichlet boundary 
condition at $t=0$,
\begin{equation}
H:=-\frac{d^2}{dt^2}+p(t)+cX'_{\omega}(t)\ ,\quad t\geq0
\label{eqn:5}
\end{equation}
can be realized as a closed symmetric operator in the Hilbert space 
$L^2(\mathbf{R}_+)$, which is self-adjoint and has purely discrete 
spectrum with probability one.
\end{thm}

Here a fractional Brownian motion $\{X_{\omega}(t)\}$ of Hurst 
parameter $h\in(0,1)$ is a stochastic process satisfying the 
following conditions:
\begin{description}
\item[(i)] $X_{\omega}(0)=0$ and $X_{\omega}(t)$ is continuous with 
respect to $t\geq0$ for all $\omega\in\Omega$;
\item[(ii)] $\{X_{\omega}(t)\}$ is a centered Gaussian process;
\item[(iii)] $\mathbf{E}[X(t)X(s)]=\frac12(t^{2h}+s^{2h}-\vert t-s\vert^{2h})$\ .
\end{description}
See Nourdin \cite{N} for details.

Since the  standard Brownian motion is the fractional Brownian motion 
of Hurst parameter $h=1/2$, Theorem 1 follows from Theorem 2 by 
letting $p(t)=t$ and $c=2/\sqrt{\beta}$.

In the following section, we give a definition of the operator $H$ as a 
closed symmetric operator. The argument there is purely deterministic. 
In \S3, we consider a random quadratic form
\begin{eqnarray}
{\cal E}_{\omega}(u,v)
&=&\int_0^{\infty}\{u'(t)\overline{v'(t)}+(p(t)+ca'_{\omega}(t))
u(t)\overline{v(t)}\}dt \nonumber\\
& &-\int_0^{\infty}c(X_{\omega}(t)-a_{\omega}(t))(u(t)\overline{v(t)})'dt\ ,
\label{eqn:6}
\end{eqnarray}
which is defined for $u,v\in L$, where
\begin{eqnarray}
L:=\{f&;&\ f(\cdot)\ \mbox{is absolutely continuous on}\ 
[0,\infty),\ f(0)=0,\ \mbox{and} \nonumber\\
& &
\int_0^{\infty}\{\vert f'(t)\vert^2+(1+p(t))\vert f(t)\vert^2\}dt<\infty\}\ .
\label{eqn:7}
\end{eqnarray}
Here we have set, as in \cite{RRV},
\begin{equation}
a_{\omega}(t)=\int_t^{t+1}X_{\omega}(s)ds\ .
\label{eqn:8}
\end{equation}
It will then be shown that for $\mathbf{P}$-almost all $\omega$, 
${\cal E}_{\omega}$ is closed, lower semi-bounded and completely 
continuous. Hence by general theory, there corresponds a 
self-adjoint operator $A_{\omega}$ which is lower semi-bounded and 
has purely discrete spectrum. To complete the proof of Theorem2, we 
then verify that this $A_{\omega}$ is a self-adjoint extension of the 
symmetric operator defined in \S2 and that this symmetric operator 
is of limit point type at $+\infty$.

\section{Definition of $H$ as a symmetric operator.}

Let $p(t)$ and $Q(t)$ be real valued continuous functions on 
$[0,\infty)$ with $Q(0)=0$, and consider the expression
\begin{equation}
H=H(p,Q)=-\frac{d^2}{dt^2}+p(t)+Q'(t)\ ,
\label{eqn:9}
\end{equation}
where $Q'(t)$ is the formal derivative of $Q(t)$.

\begin{definition}[\cite{M1}] A complex valued function $u(t)$ belongs 
to the space ${\cal C}={\cal C}(p,Q)$ if and only if $u$ is absolutely 
continuous on $[0,\infty)$ and if there are $\alpha\in\mathbf{C}$ 
and $v(\cdot)\in L^1_{loc}(\mathbf{R}_+)$ such that 
\begin{equation}
u'(t)=\alpha+Q(t)u(t)+\int_0^t
\{p(y)u(y)-Q(y)u'(y)-v(y)\}dy\ .
\label{eqn:10}
\end{equation}
\end{definition}

It is clear that $v$ is uniquely determined from $u$ as an element 
of $L^1_{loc}(\mathbf{R}_+)$. We define $Hu=H(p,Q)u=v$ for each 
$u\in{\cal C}$.

If we define, following Savchuk and Shkalikov \cite{SS}, the 
\lq\lq quasi-derivative\rq\rq of $u$ by 
$u^{[1]}(t)=u'(t)-Q(t)u(t)$, then it is easily seen that $u\in{\cal C}$ 
is equivalent to saying that $u(t)$ and $u^{[1]}(t)$ are absolutely 
continuous and satisfy
\begin{equation}
v(t):=-(u^{[1]}(t))'+(p(t)-Q^2(t))u(t)-Q(t)u^{[1]}(t)
\in L^1_{loc}(\mathbf{R}_+)\ .
\label{eqn:11}
\end{equation}
In fact this $v$ coincides with $Hu$ in Definition 1. 

The equation $Hu=\lambda u$ for $\lambda\in\mathbf{C}$ is interpreted 
as the system of integral equations
\begin{eqnarray}
u(t)&=&u(0)+\int_0^tu'(s)ds\ ; \nonumber\\
u'(t)&=&u'(0)+u(0)(Q(t)+P(t)-\lambda t) \nonumber\\
&\ &+\int_0^t\{Q(t)-Q(y)-\lambda(t-y)+P(t)-P(y)\}u'(y)dy\ ,
\label{eqn:12}
\end{eqnarray}
where $P(t):=\int_0^tp(s)ds$, 
or of differential equations
\begin{equation}
\frac{d}{dt}
\left[\begin{array}{c} u(t) \\ u^{[1]}(t) \end{array}\right]
=\left[\begin{array}{cc} Q(t) & 1 \\ 
\lambda-Q^2(t)-p(t) & -Q(t) \end{array}\right]
\left[\begin{array}{c} u(t) \\ u^{[1]}(t)\end{array}\right]\ ,
\label{eqn:13}
\end{equation}
with a given initial condition, 
both (\ref{eqn:12}) and (\ref{eqn:13}) being uniquely solvable. 
Moreover, the Green's formula is valid in the sense that for any 
$u_1,u_2\in{\cal C}$ and $0\leq a<b<\infty$,　one has
\begin{equation}
\int_a^b\{(Hu_1)(t)u_2(t)-u_1(t)(Hu_2)(t)\}dt
=[u_1,u_2](t)\vert_{t=a}^{t=b}\ ,
\label{eqn:14}
\end{equation}
where
\begin{equation}
[u_1,u_2](t):=u_1(t)u'_2(t)-u'_1(t)u_2(t)
=u_1(t)u_2^{[1]}(t)-u^{[1]}(t)u_2(t)\ .
\label{eqn:15}
\end{equation}

With these interpretations, it can be verified that $H$ can be treated 
quite similarly as the classical Sturm-Liouville operator. For example, 
if we define two spaces
\begin{equation}
{\cal D}_0:=\{u\in{\cal C}\cap L^2(\mathbf{R}_+);\ 
u(0)=0\ ,\quad Hu\in L^2(\mathbf{R}_+)\}
\label{eqn:16}
\end{equation}
and
\begin{equation}
{\cal D}_0^S:=\{u\in{\cal D}_0;\ \lim_{t\to\infty}[u,v](t)=0\ \mbox{for any}
\ v\in{\cal D}_0\}\ ,
\label{eqn:17}
\end{equation}
and if we let
\begin{equation}
H_0=H\vert_{{\cal D}_0}\ ;\quad H_0^S=H\vert_{{\cal D}_0^S}\ ,
\label{eqn:18}
\end{equation}
then we can prove that $H_0^{\ast}=H_0^S$ and 
$(H_0^S)^{\ast}=H_0$ hold in exactly the same way as the proof of 
Theorem 10.11 of \cite{St}. (From that argument, it follows in particular 
that ${\cal D}_0^S$ is dense in $L^2(\mathbf{R}_+)$.) Thus $H_0^S$ 
is a closed symmetric operator. Moreover, Weyl-Titchmarsh theory 
is valid also for $H$, and it holds that $H_0^S$ is self-adjoint if and 
only if $H$ is of limit point type at $+\infty$, which is true if and only 
if for some $\lambda\in\mathbf{C}$, the equation $Hu=\lambda u$ has a 
solution which is {\em not} square integrable near $+\infty$.

Taking $Q(t)=cX_{\omega}(t)$ in the definition of $H_0^S$ above, 
we can realize the stochastic Airy operator as a closed symmetric 
operator in $L^2(\mathbf{R}_+)$. 

We shall denote by $H_0(\omega)$ and $H_0^S(\omega)$ the operators 
defined by (\ref{eqn:18}) with this choice of $Q(t)$. Also we denote by 
$H(\omega)$ the expression (\ref{eqn:9}) itself with the same 
choice of $Q(t)$.

\begin{rem}
As is pointed out by Edelman and Sutton \cite{ES} and Bloemendal 
\cite{B}, we can transform $H$ in (\ref{eqn:9}) to a classical 
Sturm-Liouville operator
$$Lf=-\frac{d^2}{dt^2}f+pf-2Q\frac{d}{dt}f-Q^2f$$
by letting $u=f\phi$, with $\phi(t)=\exp\Bigl(\int_0^tQ(x)dx\Bigr)$. 
This makes it much clearer that Weyl-Titchmarsh theory is 
still valid for our generalized Sturm-Liouville operator $H$.

For a recent systematic treatment of generalized Sturm-Liouville 
operator including those described in \S2, see for example 
\cite{EG}.
\end{rem}

\section{Almost sure self-adjointness of the stochastic Airy operator.}

In this section, we construct a self-adjoint operator $A_{\omega}$ 
which is associated with the quadratic form ${\cal E}_{\omega}$ 
defined in (\ref{eqn:6}), and show that $A_{\omega}=H_0^S(\omega)$. 
All necessary assertions concerning ${\cal E}_{\omega}$ are 
deduced from the lemma below.

\begin{lemma}
If $\{X_{\omega}(t)\}$ is the fractional Brownian motion with Hurst 
parameter $h\in(0,1)$, then with probability one, 
$a'_{\omega}(t)={\cal O}(\sqrt{\log t})$ and 
$X_{\omega}(t)-a_{\omega}(t)={\cal O}(\sqrt{\log t})$ as $t\to\infty$, 
where $a_{\omega}(t)$ is defined by (\ref{eqn:8}).
\end{lemma}

{\bf Proof.} Let $n\leq t\leq n-1$. We then have
$$|a'_{\omega}(t)|=|X_{\omega}(t+1)-X_{\omega}(t)|
\leq|X_{\omega}(t+1)-X_{\omega}(n+1)|
+|X_{\omega}(n+1)-X_{\omega}(n)|+|X_{\omega}(t)-X_{\omega}(n)|\ ,
$$
and hence
$$\sup_{n\leq t\leq n+1}|a'_{\omega}(t)|
\leq2\sup_{0\leq s\leq1}|X_{\omega}(n+s)-X_{\omega}(n)|
+\sup_{0\leq s\leq1}|X_{\omega}(n+1+s)-X_{\omega}(n+1)|\ .
$$
Since the process $\{X_{\omega}(t+T)-X_{\omega}(T)\}_{t\geq0}$ 
has the same law as $\{X_{\omega}(t)\}_{t\geq0}$ for each $T\geq0$, 
we see that for $b>0$ and $n\geq1$,
\begin{eqnarray*}
&\ &
\mathbf{P}\Bigl(\sup_{n\leq t\leq n+1}|a'_{\omega}(t)|\geq
 b\sqrt{\log n}) \\
&\leq&
\mathbf{P}\Bigl(\sup_{0\leq s\leq1}|X_{\omega}(n+s)-X_{\omega}(n)|
\geq\frac b3\sqrt{\log n}\Bigr)
+\mathbf{P}\Bigl(\sup_{0\leq s\leq1}|X_{\omega}(n+1+s)-
X_{\omega}(n+1)|\geq\frac b3\sqrt{\log n}\Bigr) \\
&=&2\mathbf{P}\Bigl(\sup_{0\leq s\leq1}|X_{\omega}(s)|
\geq\frac b3\sqrt{\log n}\Bigr)\ .
\end{eqnarray*}
Now it is known (see Theorem 4.1 of \cite{N}) that for the fractional 
Brownian motion with Hurst parameter $h\in(0,1)$ the estimate
$$\lim_{x\to\infty}x^{-2}\log\mathbf{P}\Bigl(\sup_{0\leq s\leq1}
X_{\omega}(s)\geq x\Bigr)=-\frac12$$
holds. Hence for any $\varepsilon\in(0,1)$, we can choose a 
$K_{\varepsilon}>0$ such that 
$$\mathbf{P}\Bigl(\sup_{0\leq s\leq1}|X_{\omega}(t)|\geq x\Bigr)
\leq K_{\varepsilon}e^{-(1-\varepsilon)x^2/2}\ .$$
Hence for $b$ greater thatn $\sqrt{2/(1-\varepsilon)}$, we have 
\begin{eqnarray*}
\sum_{n\geq1}\mathbf{P}\Bigl(\sup_{n\leq t\leq n+1}|a'_{\omega}(t)|
\geq b\sqrt{\log n}\Bigr)
&\leq&
\sum_{n\geq1}2K_{\varepsilon}\exp\Bigl[-\frac{1-\varepsilon}{2}
b^2\log n\Bigr] \\
&=&2K_{\varepsilon}\sum_{n\geq1}n^{-(1-\varepsilon)b^2/2}<\infty\ ,
\end{eqnarray*}
and hence by Borel-Cantelli lemma, for $\mathbf{P}$-almost all 
$\omega\in\Omega$, we can choose a constant $C_{\omega}$ so that 
$$\sup_{n\leq t\leq n+1}|a'_{\omega}(t)|\leq C_{\omega}\sqrt{\log n}\quad 
(n\geq1)\ .$$
As to the assertion for $X_{\omega}(t)-a_{\omega}(t)$, we first 
note
$$\sup_{n\leq t\leq n+1}|X_{\omega}(t)-a_{\omega}(t)|
=\sup_{n\leq t\leq n+1}\left\vert\int_0^1
(X_{\omega}(t+s)-X_{\omega}(t))ds\right\vert
\leq\sup_{n\leq t\leq n+1}\sup_{0\leq s\leq1}
|X_{\omega}(t+s)-X_{\omega}(t)|$$
and
\begin{eqnarray*}
& &|X_{\omega}(t+s)-X_{\omega}(t)| \\
\ \\
&\leq&
\left\{
\begin{array}{rl}
|X_{\omega}(t)-X_{\omega}(n)|+|X_{\omega}(t+s)-X_{\omega}(n)| &
 (0\leq s\leq n+1-t) \\
\ \\
|X_{\omega}(t+s)-X_{\omega}(n+1)|+|X_{\omega}(n+1)-X_{\omega}(n)|+
|X_{\omega}(t)-X_{\omega}(n)| &
(n+1-t\leq s\leq1)\ .
\end{array}\right.
\end{eqnarray*}
Then as in the case of $a'_{\omega}(t)$, we again have 
$$\sup_{n\leq t\leq n+1}|X_{\omega}(t)-a_{\omega}(t)|
\leq2\sup_{0\leq s\leq1}|X_{\omega}(n+s)-X_{\omega}(n)|
+\sup_{0\leq s\leq1}|X_{\omega}(n+1+s)-X_{\omega}(n+1)|\ ,$$
and the rest of the proof is the same as the preceding argument. 

\bigskip
The following proposition summarizes all the properties of 
${\cal E}_{\omega}$ which we need.

\begin{prop} Fix an $\omega$ for which the conclusion of Lemma 1 
holds. Then under the assumption of Theorem 2, we have that
\begin{description}
\item[(i)] ${\cal E}_{\omega}(u,v)$ is well defined for all $u,v\in L$;
\item[(ii)] ${\cal E}_{\omega}$ is lower semi-bounded;
\item[(iii)] ${\cal E}_{\omega}$ is closed; 
\item[(iv)] ${\cal E}_{\omega}$ is completely continuous in the sense 
that any $L^2$-bounded sequence $\{u_n\}$ in $L$ such that 
$\sup_n{\cal E}_{\omega}(u_n,u_n)<\infty$ contains an $L^2$-convergent 
subsequence.
\end{description}
\end{prop}

{\bf Proof.} By Lemma 1, we have 
$p(t)+ca'_{\omega}(t)={\cal O}(1+p(t))$ and 
$X_{\omega}(t)-a_{\omega}(t)=o(\sqrt{1+p(t)})$ as $t\to\infty$. In 
particular for any $\varepsilon\in(0,1/c)$, one can choose a constant 
$C_1=C_1(\omega,\varepsilon)$ so that 
$|X_{\omega}(t)-a_{\omega}(t)|\leq\varepsilon\sqrt{C_1+p(t)}$ holds 
for all $t\geq0$. This together with Schwarz inequality shows that the 
integrals in (\ref{eqn:6}) are absolutely convergent for all $u,v\in L$, 
and hence ${\cal E}_{\omega}$ is well defined on $L$. To verify assertion 
(ii), we write 
\begin{eqnarray}
{\cal E}_{\omega}(u,u)
&=&\int_0^{\infty}|u'(t)|^2dt
+\int_0^{\infty}(p(t)+ca'_{\omega}(t))|u(t)|^2dt \nonumber \\
& &-c\int_0^{\infty}(X_{\omega}(t)-a_{\omega}(t))
\Bigl(u(t)\overline{u'(t)}+u'(t)\overline{u(t)}\Bigr)dt
\label{eqn:19}
\end{eqnarray}
for $u\in L$. Noting $2|xy|\leq|x|^2+|y|^2$, we see that the third 
integral on the right hand side of (\ref{eqn:19}) is bounded in 
absolute value by
$$2\varepsilon c\int_0^{\infty}\sqrt{C_1+p(t)}|u(t)||u'(t)|dt
\leq \varepsilon c\Bigl\{\int_0^{\infty}(C_1+p(t))|u(t)|^2dt
+\int_0^{\infty}|u'(t)|^2dt\Bigr\}\ ,$$
which is finite. Hence we have 
\begin{eqnarray}
{\cal E}_{\omega}(u,u)
&\geq&
(1-\varepsilon c)\int_0^{\infty}|u'(t)|^2dt
-\varepsilon cC_1\int_0^{\infty}|u'(t)|^2dt \nonumber \\
& &+\int_0^{\infty}\{(1-\varepsilon c)p(t)+ca'_{\omega}(t)\}
|u(t)|^2dt
\label{eqn:20}
\end{eqnarray} 
If we further take a $\delta\in(0,1-\varepsilon c)$, then there is a 
constant $C_2=C_2(\omega,p(\cdot),c,\varepsilon,\delta)$ such that 
\begin{equation}
\inf_{t\geq0}\{(1-\varepsilon c-\delta)p(t)+ca'_{\omega}(t)\}
\geq-C_2\ .
\label{eqn:21}
\end{equation}
Hence letting $C=\varepsilon cC_1+C_2$, we obtain the lower bound 
\begin{eqnarray}
{\cal E}_{\omega}(u,u)
&\geq&
(1-\varepsilon c)\int_0^{\infty}|u'(t)|^2dt+
\delta\int_0^{\infty}p(t)|u(t)|^2dt
-C\int_0^{\infty}|u(t)|^2dt \nonumber \\
&\geq&-C\int_0^{\infty}|u(t)|^2dt=-C(u,u)\ ,
\label{eqn:22}
\end{eqnarray}
which proves the assertion (ii).

Now fix a constant $\gamma>C$ and define the norm 
 $\Vert\cdot\Vert_{\gamma}$ on $L$ by 
$\Vert u\Vert_{\gamma}^2={\cal E}_{\omega}(u,u)+\gamma(u,u)$. 
Then by (\ref{eqn:22}),
\begin{equation}
\Vert u\Vert^2_{\gamma}\geq(1-\varepsilon c)
\int_0^{\infty}|u'(t)|^2dt+\delta\int_0^{\infty}p(t)|u(t)|^2dt
+(\gamma-C)\int_0^{\infty}|u(t)|^2dt
\label{eqn:23}
\end{equation}
holds for $u\in L$. To show that $L$ is complete under the norm 
$\Vert\cdot\Vert_{\gamma}$, let $\{u_n\}\subset L$ be a Cauchy 
sequence with respect to this norm. For this sequence, we have from 
(\ref{eqn:23})
\begin{equation}
\lim_{n,m\to\infty}\Vert u_n-u_m\Vert_2=0\ ,
\label{eqn:24}
\end{equation}
\begin{equation}
\lim_{n,m\to\infty}\Vert u'_n-u'_m\Vert_2=0,
\label{eqn:25}
\end{equation}
and
\begin{equation}
\lim_{n,m\to\infty}\int_0^{\infty}p(t)|u_n(t)-u_m(t)|^2dt=0\ .
\label{eqn:26}
\end{equation}
By (\ref{eqn:24}) and (\ref{eqn:25}), there are $u,v\in L^2(\mathbf{R}_+)$ 
such that $\Vert u_n-u\Vert_2\to0$ and $\Vert u'_n-v\Vert_2\to0$ as 
$n\to\infty$. In particular we have $\int_0^x|u'_n(t)-v(t)|dt\to0$ for 
any $x>0$. Hence we can let $n\to\infty$ in the equality 
$u_n(x)=\int_0^xu'_n(t)dt$ and we see that $\{u_n(x)\}_n$ converges 
pointwise to a continuous function which must coincide almost 
everywhere with $u(t)$, namely $u(x)=\int_0^xv(t)dt$ for almost 
every $x$. Thus $u(t)$ is absolutely continuous, $u(0)=0$, and 
$u'(t)=v(t)\in L^2(\mathbf{R}_+)$. By (\ref{eqn:26}) and Fatou's 
lemma we get 
$$\int_0^{\infty}p(t)|u(t)|^2dt\leq\liminf_{n\to\infty}
\int_0^{\infty}p(t)|u_n(t)|^2dt<\infty\ ,$$
so that $u\in L$, and 
$$\lim_{n\to\infty}\int_0^{\infty}p(t)|u_n(t)-u(t)|^2dt=0\ .$$
Since we can choose a constant $K_{\omega}$ such that
$${\cal E}_{\omega}(u,u)\leq K_{\omega}\left\{
\int_0^{\infty}|u(t)|^2dt+\int_0^{\infty}|u'(t)|^2dt
+\int_0^{\infty}p(t)|u(t)|^2dt\right\}$$
in the same way as we proved (\ref{eqn:23}), we finally obtain 
$\Vert u_n-u\Vert_{\gamma}\to0$, and hence $L$ is complete under 
$\Vert\cdot\Vert_{\gamma}$. 

If $\{u_n\}\subset L$ is a 
$\Vert\cdot\Vert_{\gamma}$-bounded sequence, then there 
are constants $M_1$, $M_2$, $M_3$ such that 
\begin{equation}
\int_0^{\infty}|u_n(t)|^2dt\leq M_1,\ \int_0^{\infty}|u_n'(t)|^2dt\leq M_2,
\ \int_0^{\infty}p(t)|u_n(t)|^2dt\leq M_3
\label{eqn:27}
\end{equation}
hold for all $n\geq1$. By the relation 
$u_n(x)^2=2\int_0^xu_n(t)u'_n(t)dt$, we have 
$$|u_n(x)|^2\leq2\int_0^{\infty}|u_n(t)||u'_n(t)|dt
\leq2\sqrt{\int_0^{\infty}|u_n(t)|^2dt}\sqrt{\int_0^{\infty}|u'_n(t)|^2dt}
\leq\sqrt{M_1M_2}\ ,$$
so that $\{u_n(x)\}$ is uniformly bounded. Since for any 
$0\leq x<y$,
$$|u_n(x)-u_n(y)|\leq\int_x^y|u'_n(t)|
\leq\sqrt{\int_x^y1dt}\sqrt{\int_0^{\infty}|u'_n(t)|^2dt}
\leq\sqrt{M_2}|y-x|$$
holds for all $n\geq1$, the sequence $\{u_n\}$ is uniformly equi-
continuous. Hence it contains a subsequence $\{u_{n'}(x)\}$ which 
converges uniformly on every compact subset of $[0,\infty)$ to a 
continuous function $u(t)$. Now for sufficiently large $X>0$ it holds 
that
$$\int_X^{\infty}|u_n(t)|^2dt\leq\frac1{\inf_{s\geq X}p(s)}\int_X^{\infty}
p(t)|u_n(t)|^2dt\leq\frac{M_3}{\inf_{s\geq X}p(s)}\ ,$$
and so
\begin{eqnarray*}
& &\int_0^{\infty}|u_{n'}(t)-u(t)|^2dt \\
&\leq&\int_0^X|u_{n'}(t)-u(t)|^2dt+\int_X^{\infty}|u(t)|^2dt
+\int_X^{\infty}|u_{n'}(t)|^2dt \\
&\leq&\int_0^X|u_{n'}(t)-u(t)|^2dt
+\liminf_{m'\to\infty}\int_X^{\infty}|u_{m'}(t)|^2dt
+\int_X^{\infty}|u_{n'}(t)|^2dt \\
&\leq&\int_0^X|u_{n'}(t)-u(t)|^2dt+\frac{2M_3}{\inf_{s\geq X}p(s)}\ .
\end{eqnarray*}
By letting $n'\to\infty$ first, and then $X\to\infty$, we obtain 
$\lim_{n'\to\infty}\Vert u_{n'}-u\Vert_2=0$.

\bigskip
We can now apply general theory, to obtain a self-adjoint operator 
$A_{\omega}$ with domain $D(A_{\omega})\subset L$ such that 
$(A_{\omega}u,v)={\cal E}_{\omega}(u,v)$ for all $u,v\in D(A_{\omega})$. 
In particular, $\psi\in D(A_{\omega})$ if and only if there is a 
$\varphi\in L^2(\mathbf{R}_+)$ such that 
${\cal E}_{\omega}(u,\psi)=(u,\varphi)$ for all 
$u\in L$ (see \cite{RS}, Theorem VIII.15). 
Moreover, the assertion (iv) above implies that $A_{\omega}$ has 
purely discrete spectrum (see \cite{RS2}, Theorem XIII.64). 
We shall show that $A_{\omega}=H_0^S(\omega)$ holds for 
every fixed $\omega$ for which the conclusion of Lemma 1 is valid.  
Let us begin by proving the following assertion.

\begin{prop}
The relation $H_0^S(\omega)\subset A_{\omega}\subset H_0(\omega)$ 
holds.
\end{prop}

{\bf Proof.} We prove $A_{\omega}\subset H_0(\omega)$ only. 
The other inclusion relation follows by taking the adjoint. 

From the definition of the quadratic form ${\cal E}_{\omega}$, it 
is easy to see that $\psi\in D(A_{\omega})$ if and only if there 
is a $\varphi\in L^2(\mathbf{R}_+)$ such that
\begin{eqnarray}
& &\int_0^{\infty}u'(t)\{\overline{\psi'(t)}-c(X_{\omega}(t)-a_{\omega}(t))
\overline{\psi(t)}\}dt \nonumber \\
&=&\int_0^{\infty}u(t)\{\overline{\varphi(t)}-(p(t)+ca'_{\omega}(t))
\overline{\psi(t)}+c(X_{\omega}(t)-a_{\omega}(t))\overline{\psi'(t)}\}dt
\label{eqn:28}
\end{eqnarray}
holds for all $u\in L$. In particular, if $\psi\in D(A_{\omega})$, 
then the relation (\ref{eqn:28}) holds for all $u\in L$ with compact 
support. Then by Lemma 2 below, we can conclude that 
$$\psi'(t)-cX_{\omega}(t)+ca_{\omega}(t)\psi(t)
=\psi^{[1]}(t)+ca_{\omega}(t)\psi(t)$$
is absolutely continuous and 
$$-(\psi^{[1]}(t)+ca_{\omega}(t)\psi(t))'
=\varphi(t)-(p(t)+ca'_{\omega}(t))\psi(t)
+c(X_{\omega}(t)-a_{\omega}(t))\psi'(t)\ ,$$
namely
\begin{equation}
\varphi(t)=-(\psi^{[1]}(t))'-cX_{\omega}(t)\psi^{[1]}(t)
+(p(t)-c^2X_{\omega}(t)^2)\psi(t)\ .
\label{eqn:29}
\end{equation}
This shows that $\psi\in{\cal D}_0(\omega)$ and 
$H_0(\omega)\psi=\varphi$.

\begin{lemma}
Let $\alpha(t)$ be locally bounded and measurable, $\gamma(t)$ be 
locally integrable on $[0,\infty)$. If 
\begin{equation}
\int_0^{\infty}u'(t)\alpha(t)dt=\int_0^{\infty}u(t)\gamma(t)dt
\label{eqn:30}
\end{equation}
holds for all $u\in L$ with compact support, then $\alpha(t)$ is 
absolutely continuous and $-\alpha'(t)=\gamma(t)$ holds almost 
everywhere on $[0,\infty)$.
\end{lemma}

{\bf Proof.} For a given $T>0$, it is easy to see that $u\in L$ has 
its support in $[0,T]$ if and only if there is a $v\in L^2(\mathbf{R}_+)$ 
such that 
\begin{equation}
u(t)=\int_0^{t\wedge T}(v(s)-\langle v\rangle_T)ds
\label{eqn31}
\end{equation}
holds, where we have set $\langle v\rangle_T=\frac1T\int_0^Tv(s)ds$\ . 
Hence (\ref{eqn:30}) becomes 
$$\int_0^T(v(t)-\langle v\rangle_T)\alpha(t)dt=
\int_0^T\left\{\int_0^t(v(s)-\langle v\rangle_T)ds\right\}
\gamma(t)dt\ ,$$
which can be rewritten in the following form
\begin{equation}
\int_0^Tv(t)(\alpha(t)-\langle\alpha\rangle_T)dt
=\int_0^Tv(t)(C(t)-\langle C\rangle_T)dt\ ,
\label{eqn:32}
\end{equation} 
where $C(t):=\int_t^T\gamma(s)ds$ . Since $v$ runs through all of 
$L^2([0,T])$, we have $\alpha(t)-\langle\alpha\rangle_T=
C(t)-\langle C\rangle_T$ almost everywhere on $[0,T]$, so that 
$\alpha'(t)=-\gamma(t)$. Since $T>0$ is arbitrary, we arrive at the 
desired conclusion.

\bigskip
To the best of the author's knowledge, it was Fukushima and Nakao 
\cite{FN} who first treated one-dimensional Schr\"odinger operator 
with singular potential including Gaussian white noise. In fact, they 
considered the operator $H_h:=-d^2/dt^2+h'(t)$ on a finite interval 
$[a,b]$ under the Dirichlet boundary condition $u(a)=u(b)=0$, 
where $h(t)$ is a bounded Borel function on $[a,b]$. They defined it 
as the self-adjoint operator $A_h$ associated to the quadratic 
form
\begin{equation}
{\cal E}_h(u,v)=\int_a^bu'(t)\overline{v(t)}dt
-\int_a^bh(t)\{u(t)\overline{v(t)}\}'dt\ ,
\label{eqn:33}
\end{equation}
which is defined on
\begin{equation}
H_0^1(a,b)=\{u;\ u(\cdot)\ \mbox{is absolutely continuous on}\ [a,b],\ 
u'(\cdot)\in L^2([a,b])\ ,\ \mbox{and}\ u(a)=u(b)=0\}
\label{eqn:34}
\end{equation}
and is closed, lower semi-bounded and completely continuous. In a 
similar manner as Lemma 2 above, it can be verified that a function 
$\psi\in L^2([a,b])$ belongs to the domain of $A_h$ if and only if 
$\psi$ and $\psi^{[1]}:=\psi'-h\psi$ are absolutely continuous on 
$[a,b]$, $\psi(a)=\psi(b)=0$ and $\varphi:=-(\psi^{[1]})'-h\psi^{[1]}-h^2\psi\in L^2([a,b])$. And in this case, we have $\varphi=A_h\psi$. 
Thus the operator $A_h$ of Fukushima and Nakao coincides with 
the operator described in \S2. 

If we introduce Pr\"ufer variables $r_{\lambda}(t)$ and 
$\theta_{\lambda}(t)$ by
\begin{equation}
\psi(t)=r_{\lambda}(t)\sin\theta_{\lambda}(t)\ ,\quad 
\psi^{[1]}(t)=r_{\lambda}(t)\cos\theta_{\lambda}(t)\ ,
\label{35}
\end{equation}
where $\psi$ is a solution to $H_h\psi=\lambda\psi$, then in 
exactly the same way as the classical Sturm-Liouville theory, 
one can prove the following oscillation theorem.

\begin{prop}
$\theta_{\lambda}(t)$ can be defined as a jointly continuous function 
of $(t,\lambda)$, which is strictly increasing in $\lambda$ for 
each fixed $t$, and which satisfies $d\theta_{\lambda}(t)/dt=1$ 
whenever $\theta_{\lambda}(t)\equiv0$ (mod $\pi$). 
If $\lambda_1<\lambda_2<\cdots$ are the eigenvalues of $A_h$,  
then the eigenfunction $\psi_h(t)$ belonging to the $k$-th 
eigenvalue $\lambda_k$ has exactly $k-1$ zeros in $(a,b)$. 
In other words, there exist $a<t_1<\cdots<t_{k-1}<b$ such that 
$\theta_{\lambda_k}(t_j)=j\pi$ ($j=1,\ldots,k-1$).
\end{prop}

\bigskip
Let us return to the proof of Theorem 2. Fix an $\omega$ for 
which the conclusion of Lemma 1 is true, and let $\lambda_1(\omega)$ 
and $\psi_{1,\omega}(t)$ be the lowest eigenvalue 
of $A_{\omega}$ and its eigenfunction. We shall show that 
$\psi_{1,\omega}(t)$ has no zeros on $(0,\infty)$, and apply the 
following theorem of Hartman \cite{H1}, to conclude that $H(\omega)$ 
is of limit-point type at $+\infty$ and that consequently 
$H_0^S(\omega)$ is self-adjoint. We note that Hartman's theorem 
is valid also for our generalized Sturm-Liouville operator.

\begin{prop} {\rm (Part (i) of Theorem in \cite{H1})} Let 
$H=H(p,Q)$ be a generalized Sturm-Liouville operator as in (\ref{eqn:9}). 
If for some real number $\mu$, the equation $H\varphi=\mu\varphi$ 
has a solution $\psi$ with only finitely many zeros on $[0,\infty)$, 
then $H$ is of limit-point type at $+\infty$.
\end{prop}

\bigskip
Let $L_0$ be the subspace of $L$ consisting of all compactly 
supported functions belonging to $L$. Then we have 
\begin{eqnarray}
\lambda_1(\omega)
&=&\inf\{{\cal E}_{\omega}(u,u);\ u\in L\ ,\quad\Vert u\Vert_2=1\} 
\nonumber \\
&\leq&\inf\{{\cal E}_{\omega};\ u\in L_0\ ,\quad \Vert u\Vert_2=1\} 
\nonumber \\
&=&\lim_{l\to\infty}
\inf\{{\cal E}_{\omega};\ u\in L_0\ ,\quad \Vert u\Vert_2=1\ ,\quad 
{\rm supp}\ u\subset[0,l]\}\ .
\label{eqn:36}
\end{eqnarray}
If we denote by $\lambda_1^l(\omega)$ the infimum appearing on the 
right hand side of (\ref{eqn:36}), then $\lambda_1^l(\omega)$ is 
nothing but the lowest eigenvalue of $H(\omega)$ on $[0,l]$ 
considered under the Dirichlet boundary condition $u(0)=u(l)=0$. 
Hence its eigenfunction $\psi_{1,\omega}^l(t)$ has no zeros on 
$(0,l)$.

Now suppose $\psi_{1,\omega}(t)$ has zeros on $(0,\infty)$ and let $l$ 
be the smallest among them, so that $\psi_{1,\omega}(t)\ne0$ for 
$0<t<l$ and $\psi_{1,\omega}(l)=0$. This shows that the function 
$\psi_{1,\omega}(t)$, considered on $0\leq t\leq l$, is the first 
eigenfunction of $H(\omega)$ on $[0,l]$ with the Dirichlet boundary 
condition at both endpoints, namely 
$\lambda_1(\omega)=\lambda_1^l(\omega)$. If 
$\theta_{\lambda}(t)$, 
$\lambda=\lambda_1(\omega)=\lambda_1^l(\omega)$, is the Pr\"ufer 
variable for $\psi_{1,\omega}(t)$, then we have 
$\theta_{\lambda}(0)=0$, $0<\theta_{\lambda}(t)<\pi$ for 
$0<t<l$, and $\theta_{\lambda}(l)=\pi$. On the other hand, 
if we take $l'>l$ and if we let $\theta_{\lambda'}(t)$, 
$\lambda':=\lambda_1^{l'}(\omega)$, be the Pr\"ufer variable for the 
first eigenfunction $\psi_{1,\omega}^{l'}(t)$ of 
$H(\omega)$ on $[0,l']$ with the Dirichlet 
boundary condition at both endpoints, then we have 
$\theta_{\lambda'}(0)=0$, $0<\theta_{\lambda'}(t)<\pi$ for 
$0<t<l'$, and $\theta_{\lambda'}(l')=\pi$. Comparing this behavior 
of $\theta_{\lambda'}(t)$ to 
that of $\theta_{\lambda}(t)$ above, we must conclude 
$\lambda'<\lambda$ from the strict monotonicity of 
$\theta_{\lambda}(t)$ with respect to $\lambda$, namely 
$\lambda_1^{l'}(\omega)<\lambda_1^l(\omega)=\lambda_1(\omega)$ 
for any $l'>l$, contradicting 
 (\ref{eqn:36}). The proof of Theorem 2 is 
now complete. 

\begin{rem}
(i) For $k\geq2$, let $\lambda_k(\omega)$ and $\psi_{k,\omega}(t)$ 
be the $k$-th eigenvalue and eigenfunction of 
$A_{\omega}=H_0^S(\omega)$, where the conclusion of Lemma 1 is valid 
for $\omega$. Then $\psi_{k,\omega}(t)$ has exactly $k-1$ zeros 
on $(0,\infty)$. This can be proved in the same way as above if we 
note the following min-max characterization of $\lambda_k(\omega)$:
\begin{equation}
\lambda_k(\omega)=\inf_M\sup\{{\cal E}_{\omega}(u,u);\ u\in M,\ 
\Vert u\Vert_2=1\}\ ,
\label{eqn:37}
\end{equation}
where $M$ runs through all $k$-dimensional linear subspaces of $L$.

\noindent
(ii) Once we have verified that $\psi_{k,\omega}(t)$ has only finitely 
many zeros for all $k$, it is now obvious that for any real $\lambda$, 
all non-trivial solution of $H(\omega)u=\lambda u$ has only finitely 
many zeros on $(0,\infty)$.

Thus from Corollary in \cite{H1}, we see that $H(\omega)$, considered 
under the boundary condition 
$u(0)\cos\theta+u'(0)\sin\theta=0$, also defines a semi-bounded 
self-adjoint operator with purely discrete spectrum. The $k$-th 
eigenfunction of this operator has exactly $k-1$ zeros on 
$(0,\infty)$.
\end{rem}

\section{Remarks and applications.}
\subsection{Explosion of the solution of a stochastic Riccati equation.}
Let $H$ be a generalized Sturm-Liouville operator described in \S2, 
and let $u(t)$ be a non-trivial solution of the equation 
$Hu=\lambda u$. If we let $z_{\lambda}(t)=u'(t)/u(t)$,  
then it is easily seen that it satisfies
\begin{equation}
z_{\lambda}(t)-z_{\lambda}(s)=Q(t)-Q(s)
+\int_s^t\{p(y)-\lambda-z_{\lambda}(y)^2\}dy
\label{eqn:38}
\end{equation}
on every interval on which $u(t)$ never vanishes, and that 
$z_{\lambda}(\tau\pm0)=\pm\infty$ whenever $u(\tau)=0$. 

If $Q(t)=cX_{\omega}(t)$, $X_{\omega}(t)$ being a sample function of the 
fractional Brownian motion of Hurst parameter $h\in(0,1)$, then 
(\ref{eqn:38}) defines a stochastic process $\{z_{\lambda,\omega}(t)\}$ 
with state space $[-\infty,+\infty]$. If, in addition, $p(t)$ satisfies 
the condition of Theorem 2, then from Remark 2, we obtain the 
following result.

\begin{prop}
For $\mathbf{P}$-almost all $\omega\in\Omega$, it holds that 
$\{z_{\lambda,\omega}(t)\}$ explodes to $-\infty$ finitely often for 
any real $\lambda$ and for any initial value $z(0)$ including 
$+\infty$.
\end{prop}

\subsection{Schr\"odinger operator with Gaussian white noise potential 
under a uniform electric field.}

In \cite{M2}, the present author considered a random Schr\"odinger 
operator on $[0,\infty)$ of the form
\begin{equation}
H_{\kappa,F}^{\theta}(\omega)
=-\frac{d^2}{dt^2}+\kappa B_{\omega}'(t)-Ft\ ,
\label{eqn:39}
\end{equation}
under the boundary condition $u(0)\cos\theta+u'(0)\sin\theta=0$ at 
the origin. Here $\{B_{\omega}(t)\}$ is the standard Brownian motion 
as before and the constants $\kappa$ and $F$ are strictly positive. 
The following result was obtained.

\begin{thm} {\rm (\cite{M2})} 
For any fixed value of $F>0$, $\kappa>0$ and $\theta\in[0,\pi)$, 
the operator 
$H_{\kappa,F}^{\theta}(\omega)$ is, with probability one, self-adjoint, 
its spectrum being $(-\infty,\infty)$. When $0<F<\kappa^2/2$, 
the spectrum of $H_{\kappa,F}^{\theta}(\omega)$ 
is almost surely pure point, while when $F\geq\kappa^2/2$, 
it is purely singular continuous.
\end{thm}

Combining Theorem 3, together with its proof in \cite{M2}, and 
Theorem 2 of the present work, we obtain the following result.

\begin{thm}
The random Schr\"odinger operator 
$$H_{\kappa,F}(\omega)=-\frac{d^2}{dt^2}+\kappa B'_{\omega}(t)
-Ft$$
on the whole space $\mathbf{R}$ is almost surely a self-adjoint 
operator in the Hilbert space $L^2(\mathbf{R})$, and its spectrum 
as a set is $(-\infty,\infty)$. $H_{\kappa,F}(\omega)$ has pure 
point spectrum with probability one when $0<F<\kappa^2/2$, 
while it has purely singular continuous spectrum if 
$F\geq\kappa^2/2$. 
\end{thm}

{\bf Proof.} Theorems 2 and 3 tells us that $H_{\kappa,F}(\omega)$ 
is almost surely of limit point type both at $+\infty$ and $-\infty$, so 
that it defines a unique self-adjoint operator almost surely. 
We know that Condition K in \cite{M2} holds for $H_{\kappa,F}(\omega)$. 
Moreover it was proved in \cite{M2} that when $0<F<\kappa^2/2$, 
the condition of Theorem 4 (ii) there is satisfied, while when 
$F\geq\kappa^2/2$, that of Theorem 4 (vi) holds. On the other hand, 
by Theorem 2, Remark 2 and Lemma in \cite{H1}, we see that with 
probability one, the equation $H_{\kappa,F}(\omega)u=\lambda u$ 
has a non-trivial solution which is square integrable near $-\infty$ 
for every real $\lambda$. Hence accordingly as $0<F<\kappa^2/2$ or 
$F\geq\kappa^2/2$, $H_{\kappa,F}(\omega)$ satisfies the condition 
(i) or (v) of Theorem 4 of \cite{M2}, and the desired conclusion 
follows.

\bigskip
{\bf Acknowledgment.} The author is grateful to Professor F. Nakano 
for bringing to his attention Bloemendal's thesis. This work was 
supported by Grants-in-Aid for Scientific Research (KAKENHI) no. 
22540205.

\end{document}